\pdfoutput=1
\documentclass[11pt,a4paper]{article}

\usepackage[T1]{fontenc}
\usepackage[utf8]{inputenc}
\usepackage{lmodern}
\usepackage{microtype}

\usepackage{geometry}
\usepackage{amsmath,amssymb}
\usepackage{graphicx}
\usepackage{float}
\usepackage{xcolor}
\usepackage{url}
\usepackage{hyperref}
\usepackage[
backend=biber,
style=numeric,
sorting=nty
]{biblatex}

\hypersetup{
	colorlinks=true,
	linkcolor=black,
	citecolor=blue!60!black,
	urlcolor=blue!60!black,
	pdftitle={Comparing Scalar Objective Functions for Multi-Criteria Engineering Optimization},
	pdfauthor={Olaf Frommann},
	pdfsubject={Comparison of scalar objective-function formulations for multi-criteria engineering optimization},
	pdfkeywords={multi-criteria optimization, scalarization, objective functions, desirability functions, fuzzy logic}
}

\title{
	Comparing Scalar Objective Functions for Multi-Criteria Engineering Optimization
}

\author{
	Olaf Frommann\\
	Institute for Aerospace Technology\\
	Hochschule Bremen, University of Applied Sciences\\
	Bremen, Germany\\
	Corresponding author: \texttt{olaf.frommann@hs-bremen.de}\\
	ORCID: \href{https://orcid.org/0009-0004-2614-4392}{0009-0004-2614-4392}
}

\date{}

\begin{document}
	
	\maketitle
	
	% ============================================================
	% Abstract
	% ============================================================
	\begin{abstract}
		Scalar objective functions are required when a multi-criteria optimization
		problem must yield a single preferred design rather than only a Pareto set. The
		choice of scalarization influences which compromise is selected, how preference
		parameters are interpreted, and whether non-supported Pareto regions can be
		reached. This paper compares four formulations for normalized bi-criteria
		minimization: weighted sums, achievement scalarizing functions, desirability
		functions, and a fuzzy-logic-based formulation. Two analytically defined Pareto
		fronts, one convex and one concave, isolate the effect of the objective
		formulation from numerical optimizer behavior. The comparison focuses on 
		reachable Pareto regions, parameter-induced selection density, compensation 
		between criteria, sensitivity, and interpretability. Results show that weighted 
		sums are simple but structurally limited on concave fronts, while achievement, 
		desirability, and fuzzy formulations reach interior non-supported regions through 
		different mechanisms. Desirability functions introduce nonlinear single-criterion
		preference mappings, whereas fuzzy rules express nonseparable and reference-dependent 
		engineering preferences.
	\end{abstract}
	
	
	% ============================================================
	\section{Introduction}
	% ============================================================
	
	Many engineering optimization problems involve competing criteria such as
	structural mass, aerodynamic drag, stability margins, manufacturability,
	robustness, and cost. The result is therefore determined not only by the
	physical model, but also by how the criteria are combined into a single
	preference statement.
	
	A Pareto front describes possible compromises between conflicting 
	criteria~\cite{Miettinen1999,Ehrgott2005}. It
	is valuable because it identifies solutions for which no criterion can be
	improved without deteriorating at least one other criterion. However, the Pareto
	front does not by itself select a single design. Many multi-criteria
	optimization methods are designed to generate or approximate such a front, but
	the actual design decision is then shifted to a subsequent selection step.
	Whenever a single design has to be selected, an explicit scalar objective
	function or an equivalent preference model is required. This becomes even more
	important when the number of criteria exceeds two, because the Pareto front
	then becomes a high-dimensional surface that is difficult to interpret.
	
	The weighted sum is one of the most common scalar objective-function
	formulations~\cite{MarlerArora2010}. Its popularity is due to its simplicity, 
	smoothness, and direct implementation. However, the practical interpretation of 
	weights can be difficult. Weight factors are often expected to express engineering
	preferences, although their effect depends strongly on criterion scaling,
	Pareto-front geometry, and the local structure of the feasible set. In
	particular, linear scalarizations cannot recover non-supported regions of a
	non-convex Pareto front~\cite{DasDennis1997,MarlerArora2010}. Nonlinear penalty 
	terms can modify this behavior, but they introduce additional shape and scaling 
	parameters and are therefore not considered here.
	
	Several alternatives have been developed to address these limitations.
	Achievement scalarizing functions and weighted Tchebycheff formulations use
	reference or ideal points and can reach Pareto-optimal points that are not
	supported by linear scalarization~\cite{Wierzbicki1986,Miettinen1999}.
	Desirability functions map each criterion to an individual preference scale 
	before aggregation. Fuzzy-logic-based formulations similarly use graded 
	preference descriptions, but additionally allow explicit rule-based interactions 
	between criteria~\cite{Zadeh1965,BellmanZadeh1970,Zimmermann1978}. Recent 
	surveys show that fuzzy multiobjective programming remains an active research 
	area in engineering and decision-support applications~\cite{KarimiEtAl2022}.
	
	This paper compares these objective-function formulations in a controlled
	setting using two analytically defined Pareto-front test cases. The simple
	bi-criteria setting is chosen to isolate the effect of the scalar objective
	function from optimizer-specific behavior, discretization error, and numerical
	noise. Although practical engineering problems often involve more criteria,
	unknown fronts, constraints, and expensive simulations, these complications
	make scalar preference modeling more difficult rather than less. The analytical
	setting therefore serves as a diagnostic case for understanding the assumptions
	and consequences of different formulations.
	
	The paper provides a diagnostic comparison of scalar objective-function
	behavior, with particular emphasis on how fuzzy rule bases change the induced
	preference mapping. The contributions are:
	\begin{itemize}
		\item a compact comparison of weighted sums, achievement scalarizing
		functions, desirability functions, and fuzzy-logic-based objective
		functions for normalized bi-criteria minimization;
		\item two analytically defined Pareto-front test cases used to distinguish
		supported and non-supported compromise regions;
		\item a reachability, selection-density, and sensitivity analysis showing how 
		selected optima move under variation of the native preference parameters of 
		each method; 
		\item a discussion of compensation, interpretability, and rule-based 
		preference modeling in engineering optimization.
	\end{itemize}
	
	% ============================================================
	\section{Scalar objective-function formulations}
	\label{sec:methods}
	% ============================================================
	
	All criteria considered in this work are normalized minimization criteria,
	
	\begin{equation}
		\mathbf{c} = (c_1,\ldots,c_n),
		\qquad
		c_i \in [0,1],
	\end{equation}
	
	where smaller values are preferable. A scalar objective-function formulation
	maps $\mathbf{c}$ to a scalar value $f(\mathbf{c})$ that is minimized. The
	methods considered below differ in how they encode preference information, how
	they allow compensation between criteria, and whether interactions between
	criteria can be represented explicitly.
	
	% ------------------------------------------------------------
	\subsection{Weighted sum}
	% ------------------------------------------------------------
	
	The weighted sum is written as
	
	\begin{equation}
		f_{\mathrm{WS}}(\mathbf{c})
		=
		\sum_{i=1}^{n} w_i c_i,
		\qquad
		w_i \ge 0,
		\qquad
		\sum_{i=1}^{n} w_i = 1.
	\end{equation}
	
	For the bi-criteria case considered in this paper, this becomes
	
	\begin{equation}
		f_{\mathrm{WS}}(c_1,c_2)
		=
		w c_1 + (1-w)c_2,
		\qquad
		w \in [0.01,0.99].
	\end{equation}
	
	The weighted sum is fully compensatory: a deterioration in one criterion can be
	offset by an improvement in another criterion if the weighted sum decreases.
	This property is often convenient, but it can be difficult to interpret in
	engineering terms. Moreover, weighted sums recover only supported Pareto points
	and therefore cannot select non-supported regions of non-convex Pareto 
	fronts~\cite{DasDennis1997,MarlerArora2010}.
	Penalty terms or nonlinear transformations can be added to weighted-sum
	objectives in order to modify the reachable set of solutions. Such extensions
	are not considered here, because they introduce additional shape and scaling
	parameters and therefore shift the focus away from the basic linear
	weighted-sum formulation.
	
	% ------------------------------------------------------------
	\subsection{Achievement scalarizing function}
	% ------------------------------------------------------------
	
	A weighted Tchebycheff formulation evaluates the largest weighted deviation
	from a reference or ideal point $\mathbf{z}$~\cite{Wierzbicki1986,Miettinen1999}. 
	In the normalized minimization
	setting used here, the ideal point is chosen as $\mathbf{z}=\mathbf{0}$. An
	augmented achievement scalarizing function is then written as
	
	\begin{equation}
		f_{\mathrm{ASF}}(\mathbf{c})
		=
		\max_i \left( w_i |c_i-z_i| \right)
		+
		\rho \sum_{i=1}^{n} w_i |c_i-z_i|,
		\qquad
		\rho > 0.
	\end{equation}
	
	For two criteria this gives
	
	\begin{equation}
		f_{\mathrm{ASF}}(c_1,c_2)
		=
		\max \left( w c_1, (1-w)c_2 \right)
		+
		\rho \left( w c_1 + (1-w)c_2 \right),
	\end{equation}
	
	with $w\in[0.01,0.99]$. In this work, $\rho=10^{-3}$ is used. The
	augmentation parameter is kept small so that the maximum weighted-deviation term
	remains dominant; its role is mainly regularization and tie-breaking. Since all
	criteria are normalized to $[0,1]$, this keeps the augmentation term small
	relative to the leading term. A systematic variation of $\rho$ would constitute
	a separate ASF tuning study and is therefore not considered here. In this
	convention, larger $w_i$ increase the penalty assigned to deviations 
	in criterion $i$. This convention should not be confused with inverse-weight 
	Tchebycheff formulations sometimes used in the multiobjective optimization 
	literature. 
	
	Compared
	with the weighted sum, this formulation is less compensatory and can recover
	Pareto-optimal points that are not supported by a linear scalarization.
	However, the weights affect the scaling of deviations from the ideal point
	rather than directly specifying qualitative engineering preferences.
	
	% ------------------------------------------------------------
	\subsection{Desirability functions}
	% ------------------------------------------------------------
	
	Desirability functions first map each criterion to an individual desirability
	value~\cite{Harrington1965,DerringerSuich1980}
	
	\begin{equation}
		d_i(c_i) \in [0,1],
	\end{equation}
	
	where $d_i=1$ represents a fully desirable value and $d_i=0$ represents a
	fully undesirable value. For a normalized minimization criterion, a simple
	family is
	
	\begin{equation}
		d_i(c_i) = 1 - c_i^{s_i},
		\qquad
		s_i > 0.
	\end{equation}
	
	The exponent $s_i$ controls the shape of the individual desirability curve.
	The case $s_i=1$ gives a linear mapping, while $s_i=2$ gives a parabolic
	mapping. In the main comparison, $s_1=s_2=2$ is used as a symmetric
	representative nonlinear configuration. This avoids criterion-specific shape
	preferences and keeps the comparison with the parabolic fuzzy memberships
	conceptually simple, although the mappings and aggregation mechanisms are not
	identical. The exponents are fixed to keep the experiment one-dimensional in
	the native sweep parameter $\alpha$; varying $s_1$, $s_2$, and $\alpha$
	jointly would constitute a separate desirability parameter study.
	
	The individual desirabilities are aggregated using a weighted geometric mean,
	as commonly used in response-surface desirability approaches~\cite{DerringerSuich1980}
	
	\begin{equation}
		D(\mathbf{c})
		=
		\left(
		\prod_{i=1}^{n} d_i(c_i)^{\alpha_i}
		\right)^{1/\sum_i \alpha_i}.
	\end{equation}
	
	For two criteria, the aggregated desirability is written as
	\begin{equation}
		D(c_1,c_2)
		=
		d_1(c_1)^{\alpha}
		d_2(c_2)^{1-\alpha},
		\qquad
		\alpha \in [0.01,0.99],
	\end{equation}
	where larger values of $D$ indicate more desirable compromises. Since all
	scalar objective functions in this paper are formulated as minimization
	objectives, the desirability formulation is represented by
	\begin{equation}
		f_{\mathrm{D}}(c_1,c_2)=1-D(c_1,c_2).
	\end{equation}
	Thus, minimizing $f_{\mathrm{D}}$ is equivalent to maximizing $D$. At the
	ideal point $c_1=c_2=0$, one obtains $D=1$ and hence
	$f_{\mathrm{D}}=0$.
	
	It should be noted that for $c_i=1$, the corresponding individual
	desirability $d_i$ becomes zero. Because the aggregation is multiplicative,
	the overall desirability $D$ also becomes zero. The resulting avoidance of
	endpoint regions is therefore not only a general feature of desirability-based
	optimization, but also a consequence of the specific single-criterion
	desirability functions and the geometric aggregation used here.
	
	Although the single-criterion desirability functions may appear intuitive, the
	method introduces two different types of parameters: shape exponents and
	aggregation weights. These parameters have different meanings and are not
	always straightforward to choose. In addition, the aggregation remains separable
	in the criteria and does not directly express rule-based interactions.
	
	% ------------------------------------------------------------
	\subsection{Fuzzy-logic-based objective function}
	% ------------------------------------------------------------
	
	The fuzzy-logic-based formulation considered here also maps each criterion to
	graded preference information, following the general concept of fuzzy sets and
	fuzzy decision-making~\cite{Zadeh1965,BellmanZadeh1970}. Fuzzy formulations have 
	also been used in multiobjective programming to encode graded preference 
	information for multiple objectives~\cite{Zimmermann1978,KarimiEtAl2022}. 
	
	For each criterion $c_i$, three membership values are evaluated,
	
	\begin{equation}
		m_{d,i}, \qquad m_{t,i}, \qquad m_{u,i},
	\end{equation}
	
	corresponding to desirable, tolerable, and undesirable criterion values. In the
	main comparison, parabolic membership functions are used because they provide a
	simple analytic representation while avoiding the derivative discontinuities of
	piecewise linear memberships. Global desirable and undesirable memberships are
	used to provide preference guidance over the full criterion interval.
	
	In contrast to ordinary desirability functions, the fuzzy formulation can also
	include explicit rules. Because the fuzzy rule base directly influences the 
	resulting scalar objective, several rule-base variants are considered. The rules 
	are chosen to distinguish between global exclusion behavior, weak reference 
	attraction, strong reference attraction, and a combination of reference attraction 
	and exclusion.
	
	For the fuzzy formulation, the native preference parameter is not a numerical
	weight. Instead, a reference compromise is defined by
	
	\begin{equation}
		c_1^* = \lambda,
		\qquad
		c_2^* = g(\lambda),
		\qquad
		\lambda \in [0.01,0.99],
	\end{equation}
	
	where $g$ denotes the analytical Pareto-front function. Thus, $\lambda$
	defines a reference point on the Pareto front. The fuzzy formulation is therefore 
	not compared as a weight-based scalarization with identical information content, 
	but as a reference-based preference model. The reference point is placed on the 
	analytical front to isolate the effect of rule activation. In practical applications, 
	such a point would typically be specified by the decision maker or updated iteratively, 
	rather than derived from a known Pareto front.
	
	The fuzzy inference yields a scalar preference value after rule aggregation and
	defuzzification. In the implementation used here, this value is affinely scaled
	to the normalized objective interval $[0,1]$, where lower values correspond to
	more preferred outcomes. The scaled output classes are represented by
	$
	q_d=0,\qquad q_t=\frac{1}{2},\qquad q_u=1,
	$
	corresponding to desirable, tolerable, and undesirable outcomes, respectively.
	The accumulated rule strengths, which include the contributions from the individual
	criterion memberships and the activated rule consequents, are denoted by 
	$R_d$, $R_t$, and $R_u$. This defines the fuzzy objective to be minimized.
	Details of the internal output membership functions, defuzzification procedure,
	and scaling are given in~\cite{Frommann2026OFMCO}.
	
	In the scaled representation used for comparison, the fuzzy objective can be
	written as
	
	\begin{equation}
		f_{\mathrm{F}}
		=
		\frac{
			q_d R_d + q_t R_t + q_u R_u
		}{
			R_d + R_t + R_u
		}.
	\end{equation}
	
	The fuzzy evaluations in this paper are computed using the open-source C++
	library \texttt{FuzzyGoal}. The present article uses one fixed fuzzy-logic 
	configuration in order to compare scalar objective-function behavior. The broader 
	fuzzy-goal framework, including linear, parabolic, and Gaussian memberships, 
	global and reference-based memberships, and min--max, soft-min--max, product, 
	and $p$-norm aggregation operators, is documented in detail in~\cite{Frommann2026OFMCO}.
	Only the specific configuration used for the numerical results is summarized here.
	
	For the present investigation, parabolic membership functions are used. The 
	desirable and undesirable functions are
	\begin{align}
		m_{d,i}(c_i)& =\left(\dfrac{c_u-c_i}{c_u-c_l}\right)^2, \\
		m_{u,i}(c_i)& =\left(\dfrac{c_i-c_l}{c_u-c_l}\right)^2,
	\end{align}
	where $c_l$ and $c_u$ denote the lower and upper criterion boundaries,
	respectively. In the normalized setting used here, $c_l=0$ and $c_u=1$. The 
	reference-based functions are defined as
	\begin{equation}
		m_{t,i}(c_i)=
		\begin{cases}
			1-\left(\dfrac{c_i^*-c_i}{c_i^*-c_l}\right)^2, & c_i \le c_i^*, \\[13pt]
			1-\left(\dfrac{c_i-c_i^*}{c_u-c_i^*}\right)^2, & c_i > c_i^*,
		\end{cases}
	\end{equation}
	where the membership takes its maximum at $c_i=c_i^*$ and vanishes exactly at the 
	interval boundaries $c_l$ and $c_u$. For the parameter sweeps used here, the 
	reference values are strictly inside the criterion interval, so that $c_l<c_i^*<c_u$.
	
	This configuration is used for all fuzzy results shown below. 
	The membership type, reference construction, and aggregation operators remain fixed. 
	The explicit rule sets F1--F4, defined in the following section, are evaluated in 
	addition to the implicit single-criterion rules of the \texttt{FuzzyGoal} library, 
	which map individual desirable, tolerable, and undesirable memberships to the 
	corresponding output classes. This ensures a well-defined fuzzy objective for all 
	evaluated criterion vectors.
		
	\subsection{Fuzzy rule-base variants}
	\label{sec:fuzzy_rule_sets}
	
	In contrast to the weighted, achievement, and desirability formulations, the
	fuzzy formulation contains an explicit rule base. The rule base is not merely
	an implementation detail; it defines how the membership values are translated
	into an overall preference statement. To illustrate this effect, four fuzzy
	rule sets are considered.
	
	The first rule set represents a global exclusion mechanism:
	\begin{equation*}
		\mathrm{F1:}\qquad
		\text{if } c_1 \text{ is undesirable or } c_2 \text{ is undesirable, then the result is undesirable.}
	\end{equation*}
	This rule penalizes solutions for which at least one criterion is classified as
	undesirable.
	
	The second rule set introduces a weak reference-dependent rule:
	\begin{equation*}
		\mathrm{F2:}\qquad
		\text{if } c_1 \text{ is tolerable and } c_2 \text{ is tolerable, then the result is tolerable.} 
	\end{equation*}
	Here, simultaneous membership in the tolerable regions does not make the
	solution desirable, but only assigns an intermediate outcome.
	
	The third rule set uses the same premise but changes the consequent:
	\begin{equation*}
		\mathrm{F3:}\qquad
		\text{if } c_1 \text{ is tolerable and } c_2 \text{ is tolerable, then the result is desirable.}
	\end{equation*}
	This rule interprets the tolerable memberships as reference-neighborhood
	memberships. If both criteria are close to their reference values, the
	corresponding compromise is treated as desirable.
	
	The fourth rule set combines the reference-attraction rule with the global
	exclusion rule:
	\begin{equation*}
		\mathrm{F4:}\qquad
		\begin{aligned}
			&\text{if } c_1 \text{ is tolerable and } c_2 \text{ is tolerable, then the result is desirable,}\\
			&\text{if } c_1 \text{ is undesirable or } c_2 \text{ is undesirable, then the result is undesirable.}
		\end{aligned}
	\end{equation*}
	This rule set combines local attraction to the reference compromise with a
	penalty for strongly undesirable individual criterion values.
	
	% ------------------------------------------------------------
	\subsection{Comparison of modeling assumptions}
	% ------------------------------------------------------------
	
	Overall, the formulations differ in the type of preference information they
	encode. Weighted sums express linear compensation, ASF formulations balance
	weighted deviations from an ideal point, desirability functions use separable
	single-criterion preference mappings with geometric aggregation, and fuzzy
	logic adds explicit rule-based interactions. The preference parameters are 
	not identical across the four methods. This is intentional: the methods 
	encode different types of preference information. The comparison therefore 
	examines how each formulation translates its native preference parameters 
	into selected Pareto compromises.
	
	% ============================================================
	\section{Analytical Pareto-front test cases}
	\label{sec:testcases}
	% ============================================================
	
	Two analytically defined bi-criteria Pareto fronts are used. Both criteria are
	normalized and minimized, and the curves represent known Pareto-optimal
	compromise sets. The setting is deliberately simple so that the effect of the
	scalar objective-function formulation can be inspected directly. In realistic
	engineering applications, additional criteria, unknown front geometry,
	constraints, and expensive or noisy simulations usually obscure these
	mechanisms. The analytical examples therefore serve as controlled diagnostic
	cases for studying how scalar objective functions select among Pareto-optimal
	compromises.
	
	% ------------------------------------------------------------
	\subsection{Convex Pareto front}
	% ------------------------------------------------------------
	
	The convex test front is
	
	\begin{equation}
		c_2 = g_{\mathrm{conv}}(c_1)
		=
		\frac{4}{3}
		\left[
		\frac{1}{(1+c_1)^2}
		-
		\frac{1}{4}
		\right],
		\qquad
		c_1\in[0,1].
		\label{eq:convex_front}
	\end{equation}
	
	This function satisfies $g_{\mathrm{conv}}(0)=1$ and
	$g_{\mathrm{conv}}(1)=0$. It is monotonically decreasing and convex on the
	considered interval.
	
	% ------------------------------------------------------------
	\subsection{Concave Pareto front}
	% ------------------------------------------------------------
	
	The concave test front is
	
	\begin{equation}
		c_2 = g_{\mathrm{conc}}(c_1)
		=
		1-(1-a)c_1-a c_1^2,
		\qquad
		c_1\in[0,1],
		\qquad
		a=0.5.
		\label{eq:concave_front}
	\end{equation}
	
	It again satisfies $g_{\mathrm{conc}}(0)=1$ and $g_{\mathrm{conc}}(1)=0$.
	The first and second derivatives are
	
	\begin{equation}
		\frac{d g_{\mathrm{conc}}}{d c_1}
		=
		-(1-a)-2a c_1,
		\qquad
		\frac{d^2 g_{\mathrm{conc}}}{d c_1^2}
		=
		-2a < 0.
	\end{equation}
	
\begin{figure}[ht]
		\centering
		\includegraphics[width=\textwidth]{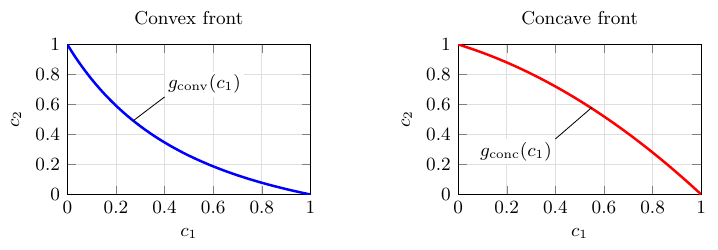}
		\caption{Analytical Pareto fronts used in this study. The convex front is
			defined by Equation~\eqref{eq:convex_front}. The concave front is defined by
			Equation~\eqref{eq:concave_front} with $a=0.5$. Both fronts are
			monotonically decreasing and represent known Pareto-optimal compromise
			sets.}
		\label{fig:fronts}
	\end{figure}

	For $0<a\le1$, the front is monotonically decreasing. Hence all points on the
	curve are non-dominated, while the interior points are non-supported by linear
	scalarization. Both Pareto fronts are shown in Fig.~\ref{fig:fronts}.
	
	% ============================================================
	\section{Numerical procedure}
	\label{sec:numerical_procedure}
	% ============================================================
	
	For each Pareto front, the curve is sampled using $10001$ uniformly distributed
	values of $c_1$. The sampling resolution was chosen such that further refinement 
	did not change the qualitative reachability and sensitivity patterns. For each 
	scalar objective formulation and each value of its native preference parameter, 
	the objective value is evaluated at all sampled points and the point with the 
	minimum objective value is recorded.
	
	The sweep parameters are
	
	\begin{equation}
		w,\alpha,\lambda \in \{0.01,0.02,\ldots,0.99\}.
	\end{equation}
	Secondary shape or regularization parameters are kept fixed: the ASF uses
	$\rho=10^{-3}$, and the desirability formulation uses $s_1=s_2=2$. This
	keeps the numerical experiment one-dimensional for each formulation and avoids
	mixing the comparison of scalarization mechanisms with a full parameter-tuning
	study of each individual method.
	
	The parameter $w$ is used for the weighted-sum and ASF formulations, whereas
	$\alpha$ is used for the desirability formulation. The parameter $\lambda$ 
	defines the fuzzy reference compromise according to Equation~\eqref{eq:convex_front} 
	or Equation~\eqref{eq:concave_front}, depending on the test case.
	
	For the main comparison with the weighted, ASF, and desirability formulations,
	rule set F3 is used:
	$
	\text{if } c_1 \text{ is tolerable and } c_2 \text{ is tolerable, then the
		overall result is desirable.}
	$
	This rule interprets the tolerable memberships as reference-neighborhood
	memberships and makes the reference compromise an active attractor in the
	scalar objective.
	
	Additional fuzzy rule sets are evaluated separately to illustrate the influence
	of the rule base on the resulting preference behavior.
	
	All data are generated with C++ programs, which are publicly available for 
	reproduction~\cite{Frommann2026Scalarization}. The fuzzy-logic-based objective
	function is evaluated using \texttt{FuzzyGoal}~\cite{Frommann2026FuzzyGoal}. The 
	resulting data files are used directly for the figures.
	
	% ============================================================
	\section{Results}
	\label{sec:results}
	% ============================================================
	
	\subsection{Reachability of Pareto-front regions}
	
\begin{figure}[ht]
		\centering
		\includegraphics[width=\textwidth]{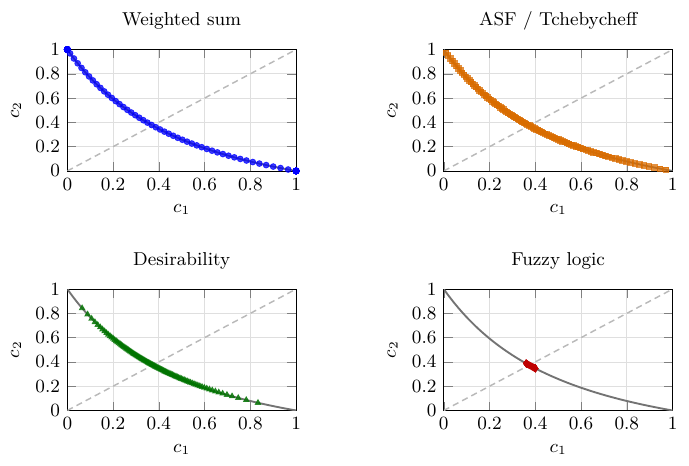}
		\caption{Selected Pareto-front points on the convex front obtained by sweeping
			the native preference parameter of each scalar objective formulation. The
			gray curve denotes the analytical Pareto front; markers denote selected
			optima. The dashed diagonal line indicates $c_1=c_2$.}
		\label{fig:reachability_convex}
	\end{figure}

\begin{figure}[H]
		\centering
		\includegraphics[width=\textwidth]{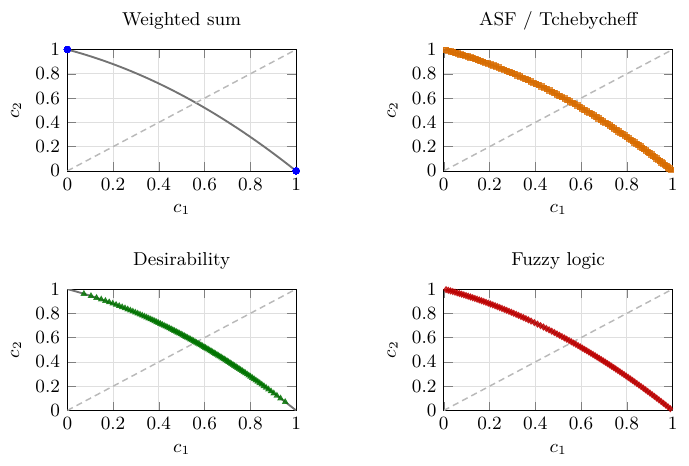}
		\caption{Selected Pareto-front points on the concave front obtained by sweeping
			the native preference parameter of each scalar objective formulation. The
			gray curve denotes the analytical Pareto front; markers denote selected
			optima. The weighted sum selects only supported points, whereas the other
			formulations can select interior regions of the concave Pareto front. The 
			dashed diagonal line indicates $c_1=c_2$.}
		\label{fig:reachability_concave}
	\end{figure}

	The Pareto-front plots in Figs.~\ref{fig:reachability_convex}
	and~\ref{fig:reachability_concave} show that the scalarizations differ not only
	in which regions of the front they can reach, but also in how densely their
	preference-parameter sweeps populate these regions. On the convex front, a
	balanced diagonal compromise exists at $c_1=c_2\approx0.373$. Several
	formulations select points near this region, but with different degrees of
	localization. On the concave front, the corresponding diagonal compromise is
	located at $c_1=c_2\approx0.562$, but the front geometry separates the
	methods more clearly: the weighted sum collapses to endpoint solutions, whereas
	ASF, desirability, and fuzzy formulations can select interior regions through
	different mechanisms.
	
\begin{figure}[ht]
		\centering
		\includegraphics[width=\textwidth]{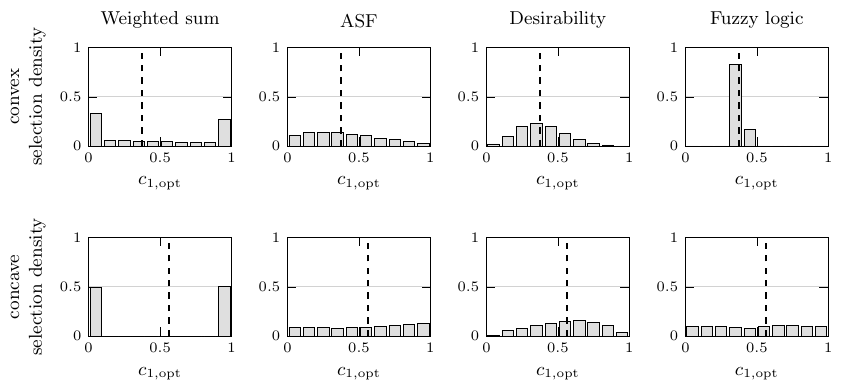}
		\caption{Parameter-induced selection density for the convex and concave Pareto
			fronts. The selected values $c_{1,\mathrm{opt}}$ are binned into ten equal
			intervals over $[0,1]$. Counts are normalized separately for each
			scalarization and front. Repeated selections of the same Pareto-front point by
			different parameter values are counted with multiplicity. The dashed vertical
			line marks the balanced compromise $c_1=c_2$.}
		\label{fig:selection_density_histograms}
	\end{figure}

	To quantify these visual clustering patterns, the selected values
	$c_{1,\mathrm{opt}}$ are binned into ten equal intervals over $[0,1]$.
	For each scalarization and each front, the bin counts are normalized by the
	number of parameter values in the corresponding sweep. Repeated selection of
	the same Pareto-front point by different parameter values is counted with
	multiplicity, because the histogram represents the parameter-induced selection
	density rather than the number of distinct Pareto points. The resulting
	histograms are shown in Fig.~\ref{fig:selection_density_histograms}.
	
	For the convex front, the weighted sum does not produce a uniform selection
	density. Instead, increased density occurs near the endpoint regions. This
	reflects the fact that weights whose supporting-line slope does not match an
	interior tangent of the convex front select endpoint optima. The ASF produces
	only a weak concentration near the balanced compromise, whereas the desirability
	formulation shows a much stronger concentration in this region. The fuzzy
	formulation is the most localized in the present configuration and selects
	points almost exclusively near the balanced compromise. All four formulations
	show a slight bias toward smaller values of $c_1$ as a result of the Pareto 
	front geometry.
	
	For the concave front, the pattern changes substantially. As expected, the
	weighted sum selects only endpoint solutions because the interior points are
	non-supported. The ASF distributes selected points over almost the entire
	front, with only a mild increase in density toward larger values of $c_1$.
	The desirability formulation avoids the endpoints and concentrates around the
	balanced compromise, with a slight bias toward larger $c_1$. The fuzzy
	formulation, in contrast to its behavior on the convex front, yields an almost
	uniform selection density over the concave front in the configuration used
	here.
	
	These results confirm that reachability and parameter-induced selection density 
	are distinct properties. A formulation may be able to reach a large part of the
	front but still assign many parameter values to specific regions. Conversely,
	two formulations may both reach interior points on a concave front while
	producing different distributions of selected optima. The selected compromise
	is therefore determined not only by whether a front region is reachable, but
	also by how the scalarization maps its preference parameters to that region.
	
\begin{figure}[ht]
		\centering
		\includegraphics[width=\textwidth]{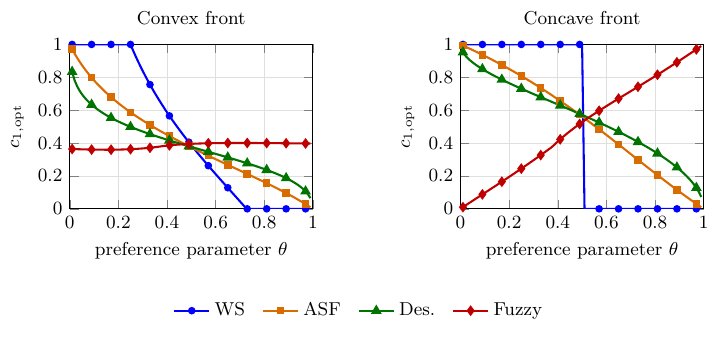}
		\caption{Sensitivity of the selected compromise to the native preference
			parameter of each formulation. The horizontal axis denotes the native
			parameter $\theta$, corresponding to $w$ for the weighted sum and ASF,
			to $\alpha$ for the desirability formulation, and to $\lambda$ for the
			fuzzy formulation.}
		\label{fig:sensitivity}
	\end{figure}

	The sensitivity curves in Fig.~\ref{fig:sensitivity} further show that the sign 
	and magnitude of the parameter response depend on the semantics of the respective 
	preference parameter. For the weighted, ASF, and desirability formulations, 
	increasing the parameter increases the emphasis on criterion $c_1$. Since both 
	criteria are minimized, this generally shifts the selected solution toward 
	smaller values of $c_1$. The response is smooth when the scalar objective selects 
	interior points continuously, but it can become discontinuous when the active 
	optimum switches between separated regions. This is most evident for the weighted 
	sum on the concave front, where the solution switches between endpoint optima.
	
	For the fuzzy formulation, the parameter $\lambda$ has a different meaning.
	It defines the location of the reference compromise directly through
	$c_1^*=\lambda$ and $c_2^*=g(\lambda)$. A positive slope of
	$c_{1,\mathrm{opt}}(\lambda)$ is therefore expected when the rule base turns
	the reference compromise into an active attractor. This behavior is visible on
	the concave front, where the fuzzy formulation selects points over a broad part
	of the front. On the convex front, the fuzzy response is much more localized
	around the balanced compromise, indicating that the reference-dependent rule is
	moderated by the global membership structure, the implicit baseline rules, and
	the geometry of the front.
	
	\subsection{Effect of the fuzzy rule base}
	\label{sec:fuzzy_rule_base_results}
	
\begin{figure}[H]
		\centering
		\includegraphics[width=\textwidth]{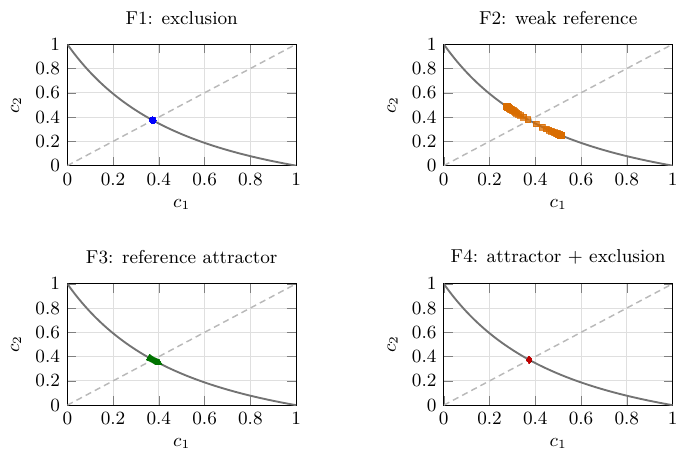}
		\caption{Effect of the fuzzy rule base on the selected Pareto-front points
			for the convex test case. The gray curve denotes the analytical Pareto
			front, while the markers show the optima obtained by sweeping the
			reference parameter $\lambda$. The dashed diagonal line indicates $c_1=c_2$.}
		\label{fig:fuzzy_rules_convex}
	\end{figure}

\begin{figure}[ht]
		\centering
		\includegraphics[width=\textwidth]{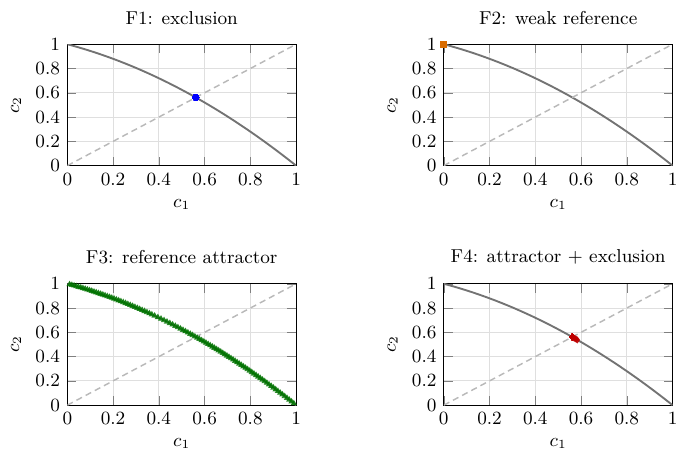}
		\caption{Effect of the fuzzy rule base on the selected Pareto-front points
			for the concave test case. The gray curve denotes the analytical Pareto front,
			while the markers show the optima obtained by sweeping the reference parameter
			$\lambda$. The dashed diagonal line indicates $c_1=c_2$.}
		\label{fig:fuzzy_rules_concave}
	\end{figure}

\begin{figure}[ht]
		\centering
		\includegraphics[width=\textwidth]{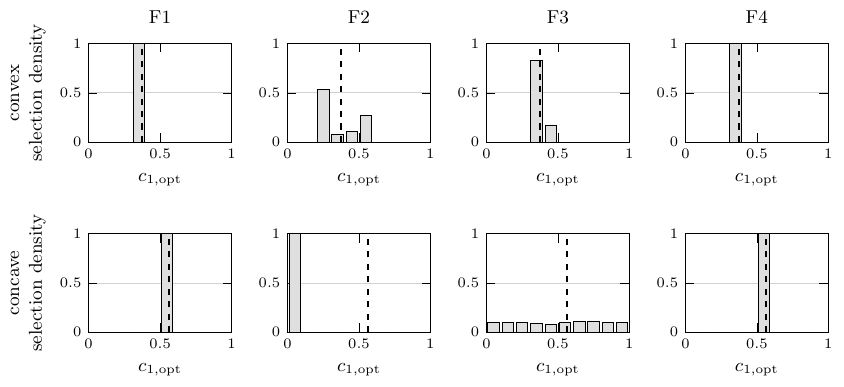}
		\caption{Parameter-induced selection density for the fuzzy rule sets F1--F4 on
			the convex and concave Pareto fronts. The selected values
			$c_{1,\mathrm{opt}}$ are binned into ten equal intervals over $[0,1]$.
			Counts are normalized separately for each rule set and front. Repeated
			selections of the same Pareto-front point by different reference values
			$\lambda$ are counted with multiplicity. The dashed vertical line marks the
			balanced compromise $c_1=c_2$.}
		\label{fig:fuzzy_rule_selection_density}
	\end{figure}

	The fuzzy formulation differs from the other scalar objective functions because
	the rule base forms an explicit part of the preference model. To illustrate
	this effect, four fuzzy rule sets were evaluated. The selected Pareto-front
	points are shown for the convex and concave fronts in
	Figs.~\ref{fig:fuzzy_rules_convex} and~\ref{fig:fuzzy_rules_concave},
	respectively. The corresponding selection-density histograms are shown in
	Fig.~\ref{fig:fuzzy_rule_selection_density}, using the same binning procedure
	as in Fig.~\ref{fig:selection_density_histograms}.
	
	The plots and histograms show that the effect of a fuzzy rule base depends
	strongly on the geometry of the Pareto front. On the convex front, all four
	rule sets select points in the vicinity of the balanced compromise
	$c_1\approx c_2$, but with different spreads. F1 and F4 produce a very narrow
	cluster around this region. F3 gives a moderately broader cluster, while F2
	produces the broadest distribution around the balanced compromise, with
	additional density toward the outer parts of the covered range. On the concave
	front, the rule sets separate more clearly. F1 and F4 again collapse to a
	narrow balanced region, F2 collapses to the endpoint at $c_1=0$, and F3
	distributes selected points almost uniformly over the front. Thus, the same
	membership functions can lead to substantially different selection densities
	depending on the rule consequents, rule interactions, and front geometry.
	
	These observations can be interpreted in terms of the mechanisms encoded by the
	individual rule sets. Rule set F1 acts mainly as a global exclusion or balancing
	mechanism. Since the explicit rule refers to undesirable memberships rather than
	to the reference-neighborhood memberships, the selected point depends only
	weakly on the reference parameter $\lambda$. The optimization therefore tends
	to avoid solutions in which one criterion becomes strongly undesirable. This
	explains the narrow clustering around the balanced compromise on both fronts.
	
	Rule set F2 contains a reference-dependent premise, but assigns only the
	intermediate tolerable consequence. On the convex front, this produces a
	partial reference-dependent response, visible as the broadest cluster among the
	four rule sets. The reference information therefore affects the selected
	optimum, but the intermediate consequent is not sufficient to turn the rule set
	into a strong reference-attraction mechanism. On the concave front, this partial
	reference effect is not robust against the front geometry, and the selected
	solutions collapse to the endpoint at $c_1=0$. Thus, a reference-dependent
	premise alone does not guarantee reference-following behavior; the rule
	consequent and the front geometry are both decisive.
	
	Rule set F3 changes the consequence of the same reference-dependent premise
	from tolerable to desirable. This turns the reference neighborhood into an
	active preference attractor. On the convex front, this broadens the selected
	region compared with F1 and F4, but the response remains centered around the
	balanced compromise. On the concave front, the effect is much stronger: the
	selected optima are distributed almost uniformly over the Pareto front. In this
	case, the fuzzy formulation can select non-supported regions when the rule base
	expresses a sufficiently strong reference-following preference.
	
	Rule set F4 combines the reference-attraction mechanism of F3 with the global
	undesirable-exclusion rule of F1. The resulting behavior is therefore more
	restrictive than that of F3 alone. In the present configuration, the exclusion
	component largely controls the selected optimum, so that F4 behaves almost like
	F1 and retains a narrow balancing tendency on both fronts. The
	reference-attraction rule modifies the local preference landscape, but does not
	by itself make the overall objective purely reference-following.
	
\begin{figure}[ht]
		\centering
		\includegraphics[width=\textwidth]{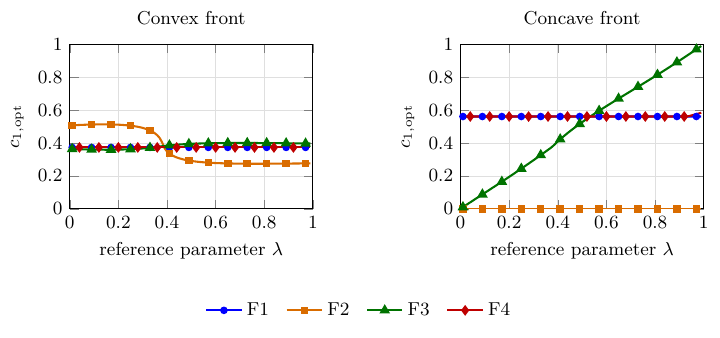}
		\caption{Sensitivity of the selected optima for the fuzzy rule sets F1--F4 to the reference parameter $\lambda$.}
		\label{fig:sensitivity_fuzzy}
	\end{figure}

	The sensitivity of the different fuzzy rule sets is shown in
	Fig.~\ref{fig:sensitivity_fuzzy}. The response depends strongly on the rule
	mechanism and on the Pareto-front geometry. F1 and F4 show little sensitivity
	to $\lambda$, because the undesirable-exclusion rule largely controls the
	selected optimum. F2 shows only partial reference dependence: on the convex
	front many parameter values are mapped to a comparatively broad region around
	the balanced compromise, whereas on the concave front the selected solutions
	collapse to $c_1=0$. In contrast, F3 exhibits an approximately linear
	dependence on $\lambda$ on the concave front, consistent with its
	reference-attraction interpretation.

	The nearly coincident sensitivity curves for F1 and F4 indicate a
	rule-dominance effect. Although F4 also contains the reference-attraction rule
	of F3, the undesirable-exclusion rule largely controls the selected optimum in
	the present configuration. Thus, F4 behaves primarily as an
	exclusion-controlled rule set rather than as a purely reference-following rule
	set. Overall, the results confirm that the preference behavior of the fuzzy
	formulation is governed by the selected rule base, by interactions between
	rules, and by the geometry of the Pareto front, rather than by a numerical
	weight.
	
	% ============================================================
	\section{Discussion}
	\label{sec:discussion}
	% ============================================================
	
	The results highlight that scalar objective-function formulations encode
	different assumptions about compensation, reachability, parameter-induced
	selection density, and preference expression. The weighted sum is simple and
	smooth, but its linear structure restricts the set of Pareto-optimal points that
	can be selected on concave fronts. Achievement scalarizing functions reduce
	this limitation by using a maximum-deviation structure, but their weights scale
	deviations from the ideal point and therefore do not directly correspond to
	simple qualitative engineering preferences. Desirability functions introduce
	nonlinear single-criterion preference mappings, which can make the treatment of
	individual criteria more intuitive, but they also add shape parameters and
	aggregation weights.
	
	The concave front illustrates that reachability alone is not sufficient to
	characterize a scalarization. The weighted sum fails to recover interior points
	because they are non-supported. The desirability formulation, by contrast,
	strongly penalizes endpoint compromises through the geometric aggregation,
	since one individual desirability becomes zero at each endpoint. The ASF and
	fuzzy formulations can reach interior points, but the distribution of selected
	points is governed by deviation balancing and rule activation, respectively.
	The histogram analysis reinforces this distinction between reachability and
	parameter-induced selection density. On the convex front, weighted sum, ASF,
	and desirability can all reach interior points, but they concentrate their
	parameter ranges in different regions. On the concave front, ASF and fuzzy
	formulations both reach non-supported regions, but they do not necessarily
	induce the same density pattern. Thus, the induced mapping from preference
	parameters to Pareto-front regions is as important as the reachable set itself.
	
	The comparison also shows that superficially similar nonlinear
	single-criterion mappings can lead to different scalar objectives. Even when
	both desirability and fuzzy formulations use parabolic single-criterion
	mappings, their behavior differs because the aggregation mechanisms are
	different. The desirability formulation remains multiplicative and separable,
	whereas the fuzzy formulation combines memberships through rule-based
	classification and aggregation. Consequently, the fuzzy objective cannot be
	characterized by its membership functions alone.
	
	The fuzzy-rule analysis demonstrates how rule consequents and rule interactions
	affect the resulting preference model. A reference-dependent premise may act
	only as a weak intermediate classification or as a strong attractor, depending
	on its consequent. In addition, rules may interact or dominate one another. The
	similarity between F1 and F4 illustrates this effect: although F4 includes the
	reference-attraction rule of F3, the undesirable-exclusion rule largely
	controls the selected optimum in the present configuration. Thus, adding a
	reference-attraction rule does not necessarily make the objective
	reference-following. The comparison between the convex and concave fronts
	further shows that rule-base effects are mediated by front geometry: a rule set
	that produces a broad balanced cluster on the convex front may collapse to an
	endpoint or become broadly reference-following on the concave front.
	
	The balanced clusters observed for several fuzzy configurations should
	therefore be interpreted as a consequence of the symmetric configuration used
	here, not as an inherent limitation of fuzzy logic. With equal criterion
	scaling, global desirable and undesirable memberships, and implicit
	single-criterion baseline rules, the objective favors regions where no
	criterion is strongly undesirable. This is meaningful when the engineering
	preference is to avoid pronounced weakness in any individual criterion. Other
	Pareto-front regions can be targeted only if the corresponding preference is
	encoded explicitly, for example through criterion-specific membership
	boundaries, reference-based desirable or undesirable memberships, stronger
	reference-attraction rules, modified exclusion rules, or additional priority
	rules.
	
	Overall, the results emphasize that scalarization should not be treated as a
	neutral numerical post-processing step. Each formulation imposes a specific
	preference structure on the Pareto front. For engineering optimization, the
	choice of scalar objective function should therefore be made together with the
	intended decision logic: linear compensation, deviation balancing, separable
	desirability aggregation, or explicit rule-based preference modeling.
	
	% ============================================================
	\section{Conclusions}
	\label{sec:conclusions}
	% ============================================================
	
	This paper compared four scalar objective-function formulations for
	multi-criteria engineering optimization using analytically controlled
	bi-criteria Pareto-front test cases. The comparison shows that no formulation
	is universally preferable. Weighted sums are simple and useful for well-behaved
	trade-offs, but are structurally limited on concave fronts. Achievement
	scalarizing functions improve reachability of non-supported regions but retain
	nontrivial weight interpretation. Desirability functions provide nonlinear
	single-criterion preference mappings but introduce additional shape and
	aggregation parameters. Fuzzy-logic-based formulations allow qualitative
	preference classes and explicit rule-based interactions between criteria.
	
	A central result is that reachability and parameter-induced selection density
	are distinct properties. A scalarization may be able to reach a Pareto-front
	region while assigning only few preference-parameter values to it, or it may
	concentrate many parameter values in a narrow region of the front. The histogram
	analysis makes this distinction explicit and shows that the same formulation
	can induce different selection-density patterns on convex and concave fronts.
	
	The fuzzy rule-base study further shows that fuzzy formulations cannot be
	characterized by their membership functions alone. Changing only the rule
	consequent can turn a reference-dependent premise from a weak intermediate
	classification into a strong preference attractor, while additional rules may
	interact or dominate one another. This makes fuzzy formulations useful when
	engineering preferences are naturally expressed through requirements,
	priorities, reference regions, or exclusion mechanisms, but it also requires
	careful rule-base design.
	
	Future work should decompose the fuzzy formulation further by separately
	examining the effects of global versus reference-based memberships, implicit
	single-criterion baseline rules, membership widths, and aggregation operators
	on the resulting preference landscape.
	
	% ============================================================
	\section*{Author contributions}
	% ============================================================
	
	Olaf Frommann: Conceptualization, Methodology, Software, Validation, Formal
	analysis, Investigation, Data curation, Writing -- original draft, Writing --
	review and editing, Visualization.
	
	% ============================================================
	\section*{Acknowledgments}
	% ============================================================
	
	The author thanks Hochschule Bremen for providing the academic environment in
	which this work was developed.
	
	% ============================================================
	\section*{Disclosure statement}
	% ============================================================
	
	The author reports there are no competing interests to declare.
	
	% ============================================================
	\section*{Funding}
	% ============================================================
	
	This research received no specific grant from any funding agency in the public,
	commercial, or not-for-profit sectors.
	
	% ============================================================
	\section*{Data availability statement}
	% ============================================================
	
	The C++ source code, generated numerical data, and PGFPlots files used to
	produce the figures are archived on Zenodo at
	\url{https://doi.org/10.5281/zenodo.20737015}~\cite{Frommann2026Scalarization}.
	The archive includes the FuzzyGoal 1.1.0 dependency snapshot used for the fuzzy
	objective-function evaluations.
	
	% ============================================================
	\section*{AI usage disclosure}
	% ============================================================
	
	Generative AI assistance was used for language editing, documentation
	structuring, drafting support, and consistency checks. The author reviewed and
	validated the final text, mathematical formulations, code, data, figures, and
	scientific claims, and remains fully responsible for the article.
	
	% ============================================================
	
	\printbibliography

@book{Ehrgott2005,
    author    = {Ehrgott, Matthias},
    title     = {Multicriteria Optimization},
    edition   = {2},
    publisher = {Springer},
    year      = {2005}
}

@book{Miettinen1999,
    author    = {Miettinen, Kaisa},
    title     = {Nonlinear Multiobjective Optimization},
    publisher = {Springer Nature},
    year      = {1998/2012},
    doi       = {10.1007/978-1-4615-5563-6}
}

@article{MarlerArora2010,
    author  = {Marler, R. Timothy and Arora, Jasbir S.},
    title   = {The Weighted Sum Method for Multi-Objective Optimization: New Insights},
    journal = {Structural and Multidisciplinary Optimization},
    volume  = {41},
    pages   = {853--862},
    year    = {2010},
    doi     = {10.1007/s00158-009-0460-7}
}

@article{DasDennis1997,
    author  = {Das, Indraneel and Dennis, John E.},
    title   = {A Closer Look at Drawbacks of Minimizing Weighted Sums of Objectives for {Pareto} Set Generation in Multicriteria Optimization Problems},
    journal = {Structural and Multidisciplinary Optimization},
    volume  = {14},
    pages   = {63--69},
    year    = {1997},
    doi     = {10.1007/BF01197559}
}

@incollection{Wierzbicki1986,
    author    = {Wierzbicki, Andrzej P.},
    title     = {On the Completeness and Constructiveness of Parametric Characterizations to Vector Optimization Problems},
    booktitle = {Multiple Criteria Decision Making Theory and Application},
    journal   = {OR Spektrum},
    volume    = {8},
    pages     = {73--87},
    year      = {1986},
    doi       = {10.1007/BF01719738}
}

@article{Harrington1965,
    author  = {Harrington, E. C.},
    title   = {The Desirability Function},
    journal = {Industrial Quality Control},
    volume  = {21},
    number  = {10},
    pages   = {494--498},
    year    = {1965}
}

@article{DerringerSuich1980,
    author  = {Derringer, George and Suich, Ronald},
    title   = {Simultaneous Optimization of Several Response Variables},
    journal = {Journal of Quality Technology},
    volume  = {12},
    number  = {4},
    pages   = {214--219},
    year    = {1980},
    doi     = {10.1080/00224065.1980.11980968}
}

@article{Zadeh1965,
    author  = {Zadeh, Lotfi A.},
    title   = {Fuzzy Sets},
    journal = {Information and Control},
    volume  = {8},
    number  = {3},
    pages   = {338--353},
    year    = {1965},
    doi     = {10.1016/S0019-9958(65)90241-X}
}

@article{BellmanZadeh1970,
    author  = {Bellman, Richard E. and Zadeh, Lotfi A.},
    title   = {Decision-Making in a Fuzzy Environment},
    journal = {Management Science},
    volume  = {17},
    number  = {4},
    pages   = {B141--B164},
    year    = {1970},
    doi     = {10.1287/mnsc.17.4.B141}
}

@article{Zimmermann1978,
    author  = {Zimmermann, H.-J.},
    title   = {Fuzzy Programming and Linear Programming with Several Objective Functions},
    journal = {Fuzzy Sets and Systems},
    volume  = {1},
    number  = {1},
    pages   = {45--55},
    year    = {1978},
    doi     = {10.1016/0165-0114(78)90031-3}
}

@article{KarimiEtAl2022,
    author  = {Karimi, N. and Feylizadeh, M. R. and Govindan, K. and Bagherpour, M.},
    title   = {Fuzzy Multi-Objective Programming: A Systematic Literature Review},
    journal = {Expert Systems with Applications},
    volume  = {192},
    pages   = {116663},
    year    = {2022},
    doi     = {10.1016/j.eswa.2022.116663}
}

@techreport{Frommann2026OFMCO,
    author      = {Frommann, Olaf},
    title       = {Objective Functions in {Multi-Criteria} Optimization:
    Weighting, {Fuzzy Logic}, and {Solution-Space} Topography},
    institution = {Hochschule Bremen},
    type        = {Technical Report},
    year        = {2026},
    note        = {Version 1.0},
    doi         = {10.5281/zenodo.20585380},
    url         = {https://doi.org/10.5281/zenodo.20585380}
}

@software{Frommann2026FuzzyGoal,
author  = {Frommann, Olaf},
title   = {{FuzzyGoal}: A {C++} Library for {Fuzzy-Logic-Based} Objective Functions},
year    = {2026},
version = {1.1.0},
doi     = {10.5281/zenodo.20593012},
url     = {https://github.com/of33/FuzzyGoal},
license = {LGPL-3.0-or-later}
}

@software{Frommann2026Scalarization,
author       = {Frommann, Olaf},
title        = {Reproducibility package for ``Comparing Scalar Objective Functions for Multi-Criteria Engineering Optimization''},
year         = {2026},
publisher    = {Zenodo},
version      = {1.0.0},
doi          = {10.5281/zenodo.20737015},
url          = {https://doi.org/10.5281/zenodo.20737015}
}
	
\end{document}